\documentclass{article}
\usepackage{amssymb,latexsym,amsfonts}

\title{Askey-Wilson relations and Leonard pairs}

\author{Raimundas Vid\= unas\thanks{Supported by the 21 Century COE Programme
"Development of Dynamic Mathematics with High Functionality" of the Ministry
of Education, Culture, Sports, Science and Technology of Japan.}\\
\em Kyushu University}
\newtheorem{theorem}{Theorem}[section]
\newtheorem{lemma}[theorem]{Lemma}

\newtheorem{definition}[theorem]{Definition}
\newtheorem{example}[theorem]{Example}
\newtheorem{remark}[theorem]{Remark}
\usepackage{latexsym}
\newcommand{\proof}{{\bf Proof. }}
\newcommand{\qed}{\hfill $\Box$}
\newcommand{\equal}{&\!\!\!=\!\!\!&}

\newcommand{\CC}{{\Bbb C}}
\newcommand{\QQ}{{\Bbb Q}}

\newcommand{\ZZ}{{\Bbb Z}}

\newcommand{\fK}{{\Bbb K}}

\date{}

\begin{document}

\maketitle

\begin{abstract} It is known that if $(A,A^*)$ is a Leonard pair, then the
linear transformations $A$, $A^*$ satisfy the Askey-Wilson relations
\begin{eqnarray*}
A^2 A^*-\beta A A^*\!A+A^*\!A^2-\gamma\left( A A^*\!+\!A^*\!A
\right)-\varrho\,A^* \equal \gamma^*\!A^2+\omega A+\eta\,I,\\
A^*{}^2\!A-\beta A^*\!AA^*\!+AA^*{}^2\!-\gamma^*\!\left(A^*\!A\!+\!A
A^*\right)-\varrho^*\!A \equal \gamma A^*{}^2\!+\omega A^*\!+\eta^*I,
\end{eqnarray*}
for some scalars %sequence of scalars
$\beta,\gamma,\gamma^*,\varrho,\varrho^*,\omega,\eta,\eta^*$. The problem of
this paper is the following: given a pair of Askey-Wilson relations as
above, how many Leonard pairs are there that satisfy those relations? It
turns out that the answer is 5 in general. We give the generic number of
Leonard pairs for each Askey-Wilson type of Askey-Wilson relations.\\

{\noindent\bf AMS 2000 MSC Classification}: 05E35, 33D45, 33C45.\\

{\noindent\bf Keywords}: Leonard pairs, Askey-Wilson relations.\\

{\noindent\bf E-mail address}: {\sf vidunas@math.kyushu-u.ac.jp}.\\
\end{abstract}

\newpage

\section{Introduction}

Throughout the paper, $\fK$ denotes an algebraically closed field. We assume
the characteristic of $\fK$ is not equal to 2. Recall that a tridiagonal
matrix is a square matrix which has nonzero entries only on the main
diagonal, on the superdiagonal and the subdiagonal. A tridiagonal matrix is
called irreducible whenever all entries on the superdiagonal and subdiagonal
are nonzero.

\begin{definition} \label{deflp} \rm
Let $V$ be a vector space over $\fK$ with finite positive dimension. By a
{\em Leonard pair} on $V$ we mean an ordered pair $(A,A^*)$, where $A:V\to
V$ and $A^*:V\to V$ are linear transformations which satisfy the following
two conditions:
\begin{enumerate}
\item[\it(i)]  There exists a basis for $V$ with respect to which the matrix
representing $A^*$ is diagonal, and the matrix representing $A$ is
irreducible tridiagonal.%
\item[\it(ii)]  There exists a basis for $V$ with respect to which the matrix
representing $A$ is diagonal, and the matrix representing $A^*$ is
irreducible tridiagonal.
\end{enumerate}
\end{definition}
\begin{remark} \rm
In this paper we do not use the conventional notation $A^*$ for the
conjugate-transpose of $A$. In a Leonard pair $(A,A^*)$, the linear
transformations $A$ and $A^*$ are arbitrary subject to the conditions {\it
(i)} and {\it (ii)} above.
\end{remark}
\begin{definition} \rm %\cite{hartwig}
Let $V,W$ be vector spaces over $\fK$ with finite positive dimensions. Let
$(A,A^*)$ denote a Leonard pair on $V$, and let $(B,B^*)$ denote a Leonard
pair on $W$. By an {\em isomorphism of Leonard pairs} we mean an isomorphism
of vector spaces $\sigma:V\mapsto W$ such that $\sigma A\sigma^{-1}=B$ and
$\sigma A^*\sigma^{-1}=B^*$. We say that $(A,A^*)$ and $(B,B^*)$ are {\em
isomorphic} if there is an isomorphism of Leonard pairs from $(A,A^*)$ to
$(B,B^*)$.
\end{definition}

Leonard pairs occur in the theory of orthogonal polynomials, combinatorics,
the representation theory of the Lie algebra $sl_2$ or the quantum group
$U_q(sl_2)$. We refer to \cite{terwgen} %\cite{TerwRome},
as a survey on Leonard pairs, and as a source of further references.

We have the following result \cite[Theorem 1.5]{TerwVid}.
\begin{theorem} \label{lptheorem}
Let $V$ denote a vector space over $\fK$ of finite positive dimension. Let
$(A,A^*)$ be a Leonard pair on $V$. Then there exists a sequence of scalars
$\beta,\gamma,\gamma^*,\varrho,\varrho^*$, $\omega,\eta,\eta^*$ taken from
$\fK$ such that
\begin{eqnarray}  \label{askwil1}
A^2 A^*-\beta A A^*\!A+A^*\!A^2-\gamma\left( A A^*\!+\!A^*\!A
\right)-\varrho\,A^* \equal \gamma^*\!A^2+\omega A+\eta\,I,\\
\label{askwil2} A^*{}^2\!A-\beta A^*\!AA^*\!+AA^*{}^2-
\gamma^*\!\left(A^*\!A\!+\!A A^*\right)-\varrho^*\!A \equal \gamma
A^*{}^2+\omega A^*\!+\eta^*I.
\end{eqnarray}
The sequence is uniquely determined by the pair $(A,A^*)$ provided the
dimension of $V$ is at least $4$.
\end{theorem}
%\vspace{6pt}
The equations (\ref{askwil1})--(\ref{askwil2}) are called the {\em
Askey-Wilson relations}. They first appeared in the work \cite{Zhidd} of
Zhedanov, where he showed that the Askey-Wilson polynomials give pairs of
infinite-dimensional matrices which satisfy the Askey-Wilson relations. We
denote this pair of equations by
$AW(\beta,\gamma,\gamma^*,\varrho,\varrho^*,\omega,\eta,\eta^*)$. We refer
to the 8 scalar parameters as the {\em Askey-Wilson coefficients}.

A natural question is the following: does a particular pair of Askey-Wilson
relations determines a Leonard pair uniquely? An example in the next section
shows that the answer is negative in general. One may ask then: if we fix
the dimension of $V$ and the 8 scalars
$\beta,\gamma,\gamma^*,\varrho,\varrho^*$, $\omega,\eta,\eta^*$, how many
Leonard pairs are there which satisfy
$AW(\beta,\gamma,\gamma^*,\varrho,\varrho^*,\omega,\eta,\eta^*)$? This is
the question that we consider in this paper.

It turns out that there may be up to 5 different Leonard pairs satisfying
the same Askey-Wilson relations. As a preliminary check, let us consider
the case $\dim V=1$. Then we have 2 equations in 2 commuting unknowns $A$
and $A^*$. Computation of a Gr\"obner basis or a resultant shows 
that there are 5 solutions in general.

Table \ref{elltab} represents our main results: the generic number of Leonard pairs,
up to isomorphism, with the same Askey-Wilson relations for various
sequences of the Askey-Wilson coefficients. 
We distinguish cases according to the classification of Askey-Wilson
relations in \cite[Section 8]{normlpaw}, which mimics Terwilliger's
classification of parameter arrays representing Leonard pairs; see
\cite{TerwLTparr} or \cite[Section 35]{terwgen} and Section \ref{lsparar}
here. These results are valid if $\dim V\ge 4$.
\begin{table}
\begin{center} \begin{tabular}{|c|c|c|}
\hline \centering Askey-Wilson coefficients & $\begin{array}{c}\mbox{Leonard}\vspace{-4pt}\\
\mbox{pairs}\end{array}$ & Askey-Wilson type\\ \hline %
$\beta\neq\pm2$, $\underline{\gamma=\gamma^*\!=0}$,
 $\widehat{\varrho}\,\widehat{\varrho}^*\!\neq 0$ & 5 & $q$-Racah \\
$\beta\neq\pm2$, $\underline{\gamma=\gamma^*\!=0}$,
 $\widehat{\varrho}=0$, $\widehat{\varrho}^*\widehat{\eta}\neq 0$ &
 4 & $q$-Hahn \\
$\beta\neq\pm2$, $\underline{\gamma=\gamma^*\!=0}$,
 $\widehat{\varrho}^*\!=0$, $\widehat{\varrho}\,\widehat{\eta}^*\!\neq 0$ &
 4 & Dual $q$-Hahn \\
$\beta\neq\pm2$, $\underline{\gamma=\gamma^*\!=0}$,
 $\widehat{\varrho}=\widehat{\eta}=0$, $\widehat{\varrho}^*\widehat{\eta}^*\neq 0$ &
 1 & $q$-Krawtchouk \\
$\beta\neq\pm2$, $\underline{\gamma=\gamma^*\!=0}$,
 $\widehat{\varrho}^*\!=\widehat{\eta}^*\!=0$, $\widehat{\varrho}\,\widehat{\eta}\neq 0$ &
 1 & Dual $q$-Krawtchouk \\
$\beta\neq\pm2$, $\underline{\gamma=\gamma^*\!=0}$,
 $\widehat{\varrho}=\widehat{\varrho}^*\!=0$, $\widehat{\eta}\,\widehat{\eta}^*\!\neq 0$ & 3 &
$\begin{array}{l} \mbox{Quantum/affine}\vspace{-5pt}\\
 \mbox{$q$-Krawtchouk}\end{array}$\\
$\beta=2$, $\gamma\,\gamma^*\neq 0$,
 $\underline{\varrho=\varrho^*=0}$ & 4 & Racah \\
$\beta=2$, $\gamma=0$, $\gamma^*\varrho\neq 0$,
 $\underline{\varrho^*=\omega=0}$ & 3 & Hahn \\
$\beta=2$, $\gamma^*\!=0$, $\gamma\,\varrho^*\!\neq 0$,
 $\underline{\varrho=\omega=0}$ & 3 & Dual Hahn \\
$\beta=2$, $\gamma=\gamma^*\!=0$, $\varrho\varrho^*\!\neq 0$,
% $\omega^2\neq\varrho\varrho^*$, %
$\underline{\eta=\eta^*\!=0}$
 & 1 & Krawtchouk \\
$\beta=-2$, $\underline{\gamma=\gamma^*=0}$,
 $\widehat{\varrho}\,\widehat{\varrho}^*\neq 0$; $\dim V$ odd & 5 & Bannai-Ito \\
$\beta=-2$, $\underline{\gamma=\gamma^*=0}$,
 $\widehat{\varrho}\,\widehat{\varrho}^*\neq 0$; $\dim V$ even & 4 & Bannai-Ito\vspace{1pt} \\
\hline
\end{tabular} \end{center}
\caption{Leonard pairs with fixed Askey-Wilson relations, if $\dim V\ge 4$}
\label{elltab}
\end{table}

The first column of Table \ref{elltab} characterizes the distinguished cases
in terms of the Askey-Wilson coefficients. The underlined expressions are
not the defining conditions; they mean that the Askey-Wilson relations can be
normalized by affine transformations
\begin{equation} \label{translation}
(A,A^*)\mapsto (tA+c, t^*\!A^*\!+c^*),\qquad \mbox{with}\
c,c^*,t,t^*\!\in\fK;\ t,t^*\neq 0,
\end{equation}
into a form where the underlined expressions hold (provided that the
preceding conditions are satisfied). Normalization of Askey-Wilson relations
is adequately explained in \cite[Section 4]{normlpaw}. %more thoroughly in Section \ref{normawrels}.
Particularly, if \mbox{$\beta\neq2$} then the Askey-Wilson relations can be
normalized so that $\gamma=0$ and \mbox{$\gamma^*\!=0$}. By
$\widehat{\varrho}$, $\widehat{\varrho}^*$, $\widehat{\omega}$,
$\widehat{\eta}$, $\widehat{\eta}^*$ we denote other Askey-Wilson
coefficients in such a normalization.

The second column indicates the generic number of Leonard pairs satisfying
Askey-Wilson relations restricted by the conditions in the first column. The
results are generic, so for some special values of the Askey-Wilson
coefficients the number of distinct Leonard pairs may be smaller. In these
special cases, one may either interpret missing Leonard pairs as degenerate, 
or one may argue that generically different Leonard pairs are isomorphic in the
special case. This is explained in Remark \ref{specialcc} and demonstrated
in Example \ref{hahnexm} here below.
If a sequence of Askey-Wilson coefficients satisfies neither condition set
of the first column, there are no Leonard pairs satisfying those
Askey-Wilson relations.

The third column gives the Askey-Wilson type of Askey-Wilson relations as
defined in \cite[Section 8]{normlpaw}. Leonard pairs have the same
Askey-Wilson type as the Askey-Wilson relations that they satisfy, according
to \cite[Theorem 8.1]{normlpaw}.

We use Terwilliger's classification of parameter
arrays representing Leonard pairs. Therefore in Section \ref{lsparar} %and Section \ref{paawrel}
we recall the definition of parameter arrays and classification terminology. In
Section \ref{normawrels} we present normalized general parameter arrays and
the Askey-Wilson relations for Leonard pairs represented by them. The
results of Table \ref{elltab} are proved in Section \ref{mainproof}.

\section{An example}
\label{almostb}

Here we give an example of Askey-Wilson relations satisfied by different
Leonard pairs. This example was observed by Curtin \cite{Curtinpr} as well.

Let $d$ be a non-negative integer, and let $V$ be a vector space with
dimension $d+1$ over $\fK$. Let $q$ denote a scalar which is not zero and
not a root of unity. Set $\beta=q^2+q^{-2}$, %q+q^{-1}$.
Notice that $\beta\neq\pm2$. 
We look for Leonard pairs on $V$ which satisfy
\begin{equation} \label{bipartite1}
AW\left(\beta,\,0,\,0,\,4-\beta^2,\,4-\beta^2,\,0,\,0,\,0\right).
\end{equation}
Existence of a Leonard pair satisfying these relations follows from
\cite{Cur2hbipT}, where Terwilliger algebras for 2-homogeneous bipartite
distance regular graphs are computed. 
The Terwilliger algebra is defined by two non-commuting generators and two
relations. The relations differ from (\ref{bipartite1}) by a scaling of the
generators. The two generators can be represented as a Leonard pair
$(A,A^*)$. The Leonard pair has the property that the tridiagonal forms for $A$ and $A^*$ 
of Definition \ref{deflp} have only zero entries on the main diagonal. A
rescaled version of $(A,A^*)$ must satisfy (\ref{bipartite1}). Besides,
Curtin \cite{Curtinpr} computed ``almost 2-homogeneous almost bipartite" Leonard
pairs satisfying the same defining relations of the Terwilliger algebra. For
these Leonard pairs, the tridiagonal forms of Definition \ref{deflp} have
precisely one nonzero entry on the main diagonal.

Here we present Leonard pairs of both kinds explicitly. They are scaled so
that they satisfy (\ref{bipartite1}). Let $A_1,A_1^*,A_2,A_2^*$ be the following matrices:
\begin{itemize}
\item $A_1$ is tridiagonal, with zero entries on the main diagonal, the
entries
\begin{equation}
\sqrt{-1}\,\frac{q^{2d-2j}-q^{2j-2d}}{q^{d-2j}+q^{2j-d}}, \qquad\mbox{for } j=0,\ldots,d-1,
\end{equation}
on the superdiagonal, and the entries
\begin{equation}
\sqrt{-1}\,\frac{q^{2j}-q^{-2j}}{q^{d-2j}+q^{2j-d}}, \qquad \mbox{for } j=1,\ldots,d,
\end{equation}
on the subdiagonal.%
\item $A_1^*$ is diagonal, with the entries
\begin{equation} %A_1^*(j,j)=
\sqrt{-1}\left(q^{d-2j}-q^{2j-d}\right),\qquad\mbox{for } j=0,\ldots,d,
\end{equation}
on the main diagonal.%
\item $A_2$ is tridiagonal, with the upper-left entry equal to
$\displaystyle \frac{q^{2d+2}-q^{-2d-2}}{q-q^{-1}}$, all other diagonal
entries equal to zero, with the entries
\begin{equation}
\frac{q^{2d-2j}-q^{2j-2d}}{q^{-2j-1}-q^{2j+1}}, \qquad\mbox{for
} j=0,\ldots,d-1,
\end{equation}
on the superdiagonal, and the entries
\begin{equation}
\frac{q^{2d+2j+2}-q^{-2d-2j-2}}{q^{2j+1}-q^{-2j-1}}, \qquad\mbox{for }
j=1,\ldots,d,
\end{equation}
on the subdiagonal.%
\item $A_2^*$ is diagonal, with the diagonal entries
$q^{2j+1}+q^{-2j-1}$, for $j=0,1,\ldots,d$.
\end{itemize}
One can routinely check that the pairs $(A_1,A_1^*)$ and $(A_2,A_2^*)$
satisfy Askey-Wilson relations (\ref{bipartite1}). Compared with the
intersection arrays for 2-homogeneous bipartite distance regular graphs in
\cite{Cur2hbipT}, we have replaced $q\mapsto q^2$, and the matrices $A_1$, $A_1^*$
are multiplied by $\sqrt{-1}\left(q^2-q^{-2}\right)\left/\left(q^{d-2}+q^{2-d}\right)\right.$.

It is a routine computation to check that the matrix pairs $(A_1,A_1^*)$ and $(A_2,A_2^*)$ satisfy
the Askey-Wilson relations (\ref{bipartite1}). Since the matrices $A_1^*$ and $A_2^*$ have different
sets of eigenvalues, the matrix pairs are not related by a conjugation. %change of base.
There are following ways to see that both $(A_1,A_1^*)$ and $(A_2,A_2^*)$ are Leonard pairs:
\begin{itemize}
\item Using Theorem 6.2 in \cite{TerwVid}. For $i=1,2$, the sufficient
conditions for $(A_i,A^*_i)$ to be a Leonard pair are the following:
\begin{itemize}
\item There exists a sequence of scalars
$\beta,\gamma,\gamma^*,\varrho,\varrho^*,\omega,\eta,\eta^*$ taken from
$\fK$ such that the Askey-Wilson relations as in {\rm (\ref{askwil1})}--{\rm (\ref{askwil2})} hold.%
\item $\widetilde{q}$ is not a root of unity, where $\widetilde{q}+\widetilde{q}^{\,-1}=\beta$.%
\item Both $A_i$ and $A^*_i$ are multiplicity free.%
\item $V$ is irreducible as an $A_i$, $A^*_i$ module.
\end{itemize}
\item By using the classification of Leonard pairs \cite[Section 35]{terwgen}.
Consider the most general $q$-Racah type:
\begin{eqnarray} \label{qracah1}
\theta_i\equal\theta_0+h\left(1-q^i\right)\left(1-s\,q^{i+1}\right)q^{-i},\\
\theta_i^*\equal\theta_0^*+h^*\left(1-q^i\right)\left(1-s^*q^{i+1}\right)q^{-i},\\
\varphi_i\equal h h^* q^{1-2i}\left(1-q^i\right)\left(1-q^{i-d-1}\right)
\left(1-r_1\,q^i\right)\left(1-r_2\,q^i\right),\\
\label{qracah9} \phi_i\equal h h^*
q^{1-2i}\left(1-q^i\right)\left(1-q^{i-d-1}\right)
\left(r_1-s^*q^i\right)\left(r_2-s^*q^i\right)\big/s^*.
\end{eqnarray}
Here $q\neq0,\pm 1$, the constants $h,h^*,s,s^*,r_1,r_2$ are nonzero and
satisfy $r_1r_2=s\,s^*q^{d+1}$, none of
$q^i,r_1q^i,r_2q^i,s^*q^i/r_1,s^*q^i/r_2$ is equal to 1 for $i=1,\ldots,d$,
and neither of $s\,q^i,s^*q^i$ is equal to 1 for $i=2,\ldots,2d$. To get the
pair $(A_1,A_1^*)$, we must replace $q\mapsto q^2$ and take
\begin{eqnarray}
\theta_0=\theta_0^*=\sqrt{-1}\!\left(q^d-q^{-d}\right),\quad
h=h^*=\sqrt{-1}\,q^d,\\
s=s^*=-q^{-2d-2},\qquad r_1=-r_2=\sqrt{-1}\,q^{-d-1},
\end{eqnarray}
and use explicit expressions in \cite[Section 27]{terwgen}. To get the pair
$(A_2,A_2^*)$, we must replace $q\mapsto q^2$ and take
\begin{eqnarray}
\theta_0=\theta_0^*=q+q^{-1},\!\quad h=h^*=q^{-1},\!\quad
s=s^*=1,\!\quad r_1=-1,\!\quad r_2=-q^{2d+2}.
\end{eqnarray}
\item By exhibiting explicit transition matrices
%from the base mentioned in part {\em (i)} of Definition \ref{deflp}
to a base mentioned in part {\em (ii)} of Definition \ref{deflp}. Entries of
the transition matrices are $q$-hypergeometric series; see \cite[Section 16]{Terw24},
\cite[Section 19]{TerwRacah} or \cite[Section 24]{terwgen}. 
Let $P_1$ denote the $(d+1)\times(d+1)$ matrix with the $(i,j)$-th entry equal to
\begin{eqnarray} \label{qracahtrz1}
(-1)^{j}\,q^{2dj}\!\left(1+q^{4j-2d}\right)\!\frac{(q^{-2d};q^2)_j(-q^{-2d};q^2)_j}
{(q^2;q^2)_j(-q^2;q^2)_j}\times\hspace{100pt}\nonumber\\ 
{}_4\phi_3\left({q^{-2i},\,q^{-2j},\,-q^{2i-2d},\,-q^{2j-2d}
\atop q^{-2d},\,\sqrt{-1}\,q^{1-d},\,-\sqrt{-1}\,q^{1-d}};\; q^2;\; q^2\;\right),
\end{eqnarray}
and let $P_2$ denote the $(d+1)\times(d+1)$ matrix with the $(i,j)$-th entry equal to
\begin{eqnarray} \label{qracahtrz2}
\frac{(1-q^{4j+2})}{q^{2(di+i-dj)}}
\frac{(-q^{2d+4};q^2)_i\,(q^{-2d};q^2)_j}
{(-q^{2-d};q^2)_i\,(q^{2d+4};q^2)_j}\,
{}_4\phi_3\!\left({q^{-2i},\,q^{-2j},\,q^{2i+2},\,q^{2j+2}
\atop q^{-2d},\,-q^{2d+4},\,-q^2};\; q^2;\; q^2\right).
\end{eqnarray}
In these expressions, $i,j\in\{0,1,\ldots,d\}$. The $q$-hypergeometric ${}_4\phi_3$
series can be written as $q$-Racah polynomials; see \cite[Section 3.2]{koekswart}.
Using $q$-difference relations for $q$-Racah polynomials, we routinely
check that $A_iP_i=P_iA_i^*$ and $A_i^*P_i=P_iA_i$ for $i=1,2$. This implies
that conjugation by $P_i$ converts the pair $(A_i,A^*_i)$ to the matrix pair
$(A_i^*,A_i)$, and condition {\em (ii)} of Definition \ref{deflp} is
satisfied. 
%\item For a dense set of values of $q$,  %
\end{itemize}
The conclusion is that $(A_1,A_1^*)$ and $(A_2,A_2^*)$ are non-isomorphic
Leonard pairs (in general), and they satisfy the same Askey-Wilson relations (\ref{bipartite1}).
Both Leonard pairs are self-dual.

Table \ref{elltab} predicts 5 Leonard pairs satisfying (\ref{bipartite1}).
Indeed, the 5 Leonard pairs are
\begin{equation} \label{bipsols}
(A_1,A_1^*),\quad (A_2,A_2^*),\quad (A_2,-A_2^*),\quad (-A_2,A_2^*),\quad
(-A_2,-A_2^*).
\end{equation}
The last 4 Leonard pairs are non-isomorphic Leonard pairs
related by affine transformations (which are affine scalings by $-1$). 
The same affine scalings of $(A_1,A^*_1)$ are isomorphic to $(A_1,A^*_1)$.
Surely, the scalings leave the Askey-Wilson relations invariant.

The complex conjugation of $\sqrt{-1}$ has the effect of multiplying both $A_1$ and $A_1^*$ by $-1$. 
The same rescaling of $(A_2,A^*_2)$ is achieved by the substitution $q\to -q$. 
The substitution $q\mapsto 1/q$
preserves the Askey-Wilson relations as well; it has the mentioned
affine rescaling action on $(A_1,A_1^*)$, and it leaves $(A_2,A^*_2)$ invariant.

\section{Leonard pairs and parameter arrays}
\label{lsparar}

Leonard pairs are represented and classified by parameter arrays. More
precisely, parameter arrays are in one-to-one correspondence with {\em
Leonard systems} \cite[Definition 3.2]{terwgen}, and to each Leonard pair
one associates 4 Leonard systems or parameter arrays. % \cite[{terwgen}
From now on, let $d$ be a non-negative integer, and let $V$ be a vector
space with dimension $d+1$ over $\fK$.
\begin{definition} \label{defpa} \rm
By a {\em parameter array} over $\fK$, of diameter $d$, we mean a sequence
\begin{equation} \label{paraar}
(\theta_0,\theta_1,\ldots,\theta_d;\;\theta_0^*,\theta_1^*,\ldots,\theta_d^*;\;
\varphi_1,\ldots,\varphi_d;\;\phi_1,\ldots,\phi_d)
\end{equation}
of scalars taken from $\fK$, that satisfy the following conditions:
\begin{enumerate}
\item[PA1.] $\theta_i\neq\theta_j$ and $\theta_i^*\neq\theta_j^*$ if $i\neq j$, for $0\le i,j\le d$.%
\item[PA2.] $\varphi_i\neq 0$ and $\phi_i\neq 0$, for $1\le i\le d$.%
\item[PA3.] $\displaystyle\varphi_i=\phi_1\sum_{j=0}^{i-1}
\frac{\theta_j-\theta_{d-j}}{\theta_0-\theta_d}+\left(\theta^*_i-\theta_0^*\right)
\left(\theta_{i-1}-\theta_d\right)$, for $1\le i\le d$.%
\item[PA4.] $\displaystyle\phi_i=\varphi_1\sum_{j=0}^{i-1}
\frac{\theta_j-\theta_{d-j}}{\theta_0-\theta_d}+\left(\theta^*_i-\theta_0^*\right)
\left(\theta_{d-i+1}\!-\theta_0\right)$, for $1\le i\le d$.%
\item[PA5.] The expressions
\[ \frac{\theta_{i-2}-\theta_{i+1}}{\theta_{i-1}-\theta_i},\qquad
\frac{\theta^*_{i-2}-\theta^*_{i+1}}{\theta^*_{i-1}-\theta^*_i}
\]
are equal and independent of $i$, for $2\le i\le d-1$.
\end{enumerate}
\end{definition}
\medskip

To get a Leonard pair from parameter array (\ref{paraar}), one must choose a
basis for $V$ and define the two linear transformations by the following
matrices (with respect to that basis):
\begin{equation} \label{split1}
\left(\begin{array}{ccccc} \theta_0 \\ 1 & \theta_1 \\ & 1 & \theta_2 \\
& & \ddots & \ddots \\ & & & 1 & \theta_d \end{array}\right),\qquad
\left(\begin{array}{ccccc} \theta^*_0 & \varphi_1 \\ & \theta^*_1 &
\varphi_2 \\ & & \theta^*_2 & \ddots \\
& & & \ddots  & \varphi_d \\ & & & & \theta^*_d \end{array}\right).
\end{equation}
Alternatively, the following two matrices define an isomorphic Leonard pair
on $V$:
\begin{equation} \label{split2}
\left(\begin{array}{ccccc} \theta_d \\ 1 & \theta_{d-1} \\ & 1 & \theta_{d-2} \\
& & \ddots & \ddots \\ & & & 1 & \theta_0 \end{array}\right),\qquad
\left(\begin{array}{ccccc} \theta^*_0 & \phi_1 \\ & \theta^*_1 &
\phi_2 \\ & & \theta^*_2 & \ddots \\
& & & \ddots  & \phi_d \\ & & & & \theta^*_d \end{array}\right).
\end{equation}

Conversely, if $(A,A^*)$ is a Leonard pair on $V$, there exists
\cite[Section 21]{terwgen} a basis for $V$ with respect to which the
matrices for $A$, $A^*$ have the bidiagonal forms in (\ref{split1}),
respectively. There exists another basis for $V$ with respect to which the
matrices for $A$, $A^*$ have the bidiagonal forms in (\ref{split2}),
respectively, with the same scalars $\theta_0,\theta_1,\ldots,\theta_d;
\theta^*_0,\theta^*_1,\ldots,\theta^*_d$. Then the following 4 sequences are
parameter arrays of diameter $d$:
\begin{eqnarray} \label{parray1}
(\theta_0,\theta_1,\ldots,\theta_d;\;\theta_0^*,\theta_1^*,\ldots,\theta_d^*;\;
\varphi_1,\ldots,\varphi_d;\;\phi_1,\ldots,\phi_d),\\ \label{parray2}
(\theta_0,\theta_1,\ldots,\theta_d;\;\theta_d^*,\ldots,\theta_{1}^*,\theta_0^*;\;
\phi_d,\ldots,\phi_1;\;\varphi_d,\ldots,\varphi_1),\\ \label{parray3}
(\theta_d,\ldots,\theta_{1},\theta_0;\;\theta_0^*,\theta_1^*,\ldots,\theta_d^*;\;
\phi_1,\ldots,\phi_d;\;\varphi_1,\ldots,\varphi_d),\\ \label{parray9}
(\theta_d,\ldots,\theta_{1},\theta_0;\;\theta_d^*,\ldots,\theta_{1}^*,\theta_0^*;\;
\varphi_d,\ldots,\varphi_1;\;\phi_d,\ldots,\phi_1).
\end{eqnarray}
If we apply to any of these parameter arrays the construction above, we get
back a Leonard pair isomorphic to $(A,A^*)$. These are all parameter arrays
which correspond to $(A,A^*)$ in this way.

The parameter arrays in (\ref{parray1})--(\ref{parray9}) are related by
permutations. The permutation group is isomorphic to $\ZZ_2\times\ZZ_2$. The
group action is without fixed points, since the eigenvalues $\theta_i$'s (or
$\theta_i^*$'s) are distinct. Let $\downarrow$ and $\Downarrow$ denote the
permutations which transform (\ref{parray1}) into, respectively,
(\ref{parray2}) and (\ref{parray3}). Observe that the composition
$\downarrow\Downarrow$ transforms (\ref{parray1}) into (\ref{parray9}). We
refer to the permutations $\downarrow$, $\Downarrow$ and
$\downarrow\Downarrow$ as {\em relation operators}, because in \cite[Section
4]{terwgen} the parameter arrays in (\ref{parray1})--(\ref{parray9})
corresponding to $(A,A^*)$ and the 4 similar parameter arrays corresponding
to the Leonard pair $(A^*,A)$ are called {\em relatives} of each other.

Parameter arrays are classified by Terwilliger in
\cite{TerwLTparr}; alternatively, see \cite[Section 35]{terwgen}.
For each parameter array, certain orthogonal polynomials
naturally occur in entries of the transformation matrix between two bases
characterized in Definition \ref{deflp} for the corresponding Leonard pair.
Terwilliger's classification largely mimics the terminating branch of
orthogonal polynomials in the Askey-Wilson scheme \cite{koekswart}.
Specifically, the classification comprises Racah, Hahn, Krawtchouk
polynomials and their $q$-versions, plus Bannai-Ito and orphan polynomials.
Classes of parameter arrays can be identified by the type of corresponding
orthogonal polynomials; we refer to them as {\em Askey-Wilson types}. The
type of a parameter array is unambiguously defined if $d\ge 3$. %$\dim V\ge 4$.
We recapitulate Terwilliger's classification in Section \ref{normawrels} by
giving general normalized parameter arrays of each type.

By inspecting Terwilliger's general parameter arrays \cite[Section
35]{terwgen}, one can observe that the relation operators $\downarrow$,
$\Downarrow$, $\downarrow\Downarrow$ do not change the Askey-Wilson type of
a parameter array (but only the free parameters such as $q,h,h^*,s$ there),
except that the $\Downarrow$ and $\downarrow\Downarrow$ relations mix up the
quantum $q$-Krawtchouk and affine $q$-Krawtchouk types. Consequently, given
a Leonard pair, all 4 associated parameter arrays have the same type, except
when parameter arrays of the quantum $q$-Krawtchouk or affine $q$-Krawtchouk
type occur. Therefore we can use the same classifying terminology for
Leonard pairs, except that we have to merge the quantum $q$-Krawtchouk and
affine $q$-Krawtchouk types. 

Expressions for Askey-Wilson coefficients in terms of parameter arrays are
given in \cite[Theorem 4.5 and Theorem 5.3]{TerwVid} and \cite[formulas
(11)-(23)]{normlpaw}. For example, we have
\begin{eqnarray} \label{betap1}
\beta+1\equal\frac{\theta_{i-2}-\theta_{i+1}}{\theta_{i-1}-\theta_i}=
\frac{\theta^*_{i-2}-\theta^*_{i+1}}{\theta^*_{i-1}-\theta^*_i},
\hspace{56pt}\mbox{for $i=2,\ldots,d-1$};\\
\label{betap2}\gamma\equal\theta_{i-1}-\beta\theta_i+\theta_{i+1},
\hspace{97pt}\mbox{for $i=1,\ldots,d-1$};\\
\label{betap3}\gamma^*\equal\theta^*_{i-1}-\beta\theta^*_i+\theta^*_{i+1},
\hspace{95pt}\mbox{for $i=1,\ldots,d-1$};\\
\label{betap4}\varrho\equal\theta_i^2-\beta\,\theta_i\,\theta_{i-1}
+\theta_{i-1}^2-\gamma\,(\theta_i+\theta_{i-1}),\hspace{18pt}\mbox{for $i=1,\ldots,d$};\\
\label{betap5}\varrho^*\equal\theta^{*\,2}_i-\beta\theta^*_i\theta^*_{i-1}
+\theta^{*2}_{i-1}- \gamma^*(\theta^*_i+\theta^*_{i-1}),\quad\mbox{for
$i=1,\ldots,d$}.
\end{eqnarray}
In principle, these equations can be used to compute parameter arrays (and
consequently, Leonard pairs) satisfying fixed Askey-Wilson relations. For
instance, one can use (\ref{betap2})--(\ref{betap3}) to eliminate
consequently $\theta_2,\ldots,\theta_d$ and $\theta_2^*,\ldots,\theta^*_d$.
Each solution of obtained equations represents a parameter array in general.
Since we are interested in counting Leonard pairs rather than parameter
arrays, we would get $5\times4=20$ solutions in general. To get an equation
system whose solutions correspond directly to Leonard pairs, one should find
$\Downarrow$-$\downarrow$-invariant equations and rewrite them in terms of
invariants of the $\Downarrow$-$\downarrow$-action. Examples of such
invariants are, for $i=0,1,\ldots,\left\lfloor\frac{d-1}2\right\rfloor$:
\[ %begin{equation} \label{z2z2invs}
\theta_i+\theta_{d-i},\quad \theta_i\theta_{d-i},\quad %\theta^*_i\theta^*_{d-i},\,\\
\varphi_i\!\left(\theta_{d-i+1}\!-\theta_{d-i}\right)+
\varphi_{d-i+1}\!\left(\theta_{i}\!-\theta_{i-1}\right),\quad
\theta^*_i\!+\theta^*_{d-i}.
\] %end{equation}
The Askey-Wilson coefficients are invariants as well. These direct equations
can be investigated and solved if $d$ is fixed and small. In general, it
seems that one cannot avoid use of explicit solutions of recurrence
relations such as (\ref{betap2})--(\ref{betap3}), which basically leads to
classification of parameter arrays. Therefore we openly use Terwilliger's
classification. In Section \ref{normawrels} we present general normalized
parameter arrays and Askey-Wilson relations for them.

\begin{remark} \rm \label{specialcc}
For non-generic instances of Askey-Wilson relations,
the number of distinct Leonard pairs is %may be %different from
smaller than the respective generic number in Table \ref{elltab}.
Within intersection theory (or moduli space) philosophy, there may be following
``reasons" for this:
\begin{itemize}
\item Some solutions of a defining equation system represent ``degenerate" objects
rather than genuine Leonard pairs.  In our
situation, degenerate objects are represented by ``parameter arrays" which do not
satisfy the conditions PA1 and PA2 of Definition \ref{defpa}.%
\item General Leonard pairs in parametrized families are supposed to be
generically different and non-isomorphic, but they may coincide or be
isomorphic for special values of the parameters, or for special instances of Askey-Wilson relations.
In these situations, one can assign a {\em multiplicity} to each
solution so that multiplicities of all solutions add up to the generic
number. Multiplicities should be defined by considering the defining
equation system locally, or by an appropriate infinitesimal deformation of the parameters.
\end{itemize}
Example \ref{hahnexm} here below presents instances of these situations.
More generally, we may expect other two standard complications:
\begin{itemize}
\item Some ``missing" solutions are at the ``infinity" (or more technically, on a
compactification of the ``moduli space" of possible Leonard pairs). We do not need this
interpretation within each Askey-Wilson type, unless we wish to have the
most generic number of 5 Leonard pairs each time.%
\item A specialized defining equation system defines an algebraic variety of
positive dimension. In this case we would have infinitely many solutions, continuous
families of them. But this situation is not actual to us. (Lemma 4.1 in \cite{normlpaw}
suggests this situation for the Askey-Wilson relations with $\beta=2$,
$\gamma=0$, $\gamma^*\!=0$, $\omega^2=\varrho\varrho^*$, but then all solutions
are degenerate if $d\ge 3$.)
\end{itemize}
\end{remark}

\section{Normalized Leonard pairs}
\label{normawrels}

Let $(A,A^*)$ denote a Leonard pair, and let $c,c^*,t,t^*$ denote scalars in
$\fK$. It is easy to see that if $t$ and $t^*$ are nonzero, then
$(t\,A+c,\,t^*A^*\!+c^*)$ is a Leonard pair again. We identify here affine
transformations (\ref{translation}) acting on Leonard pairs. A corresponding
action on parameter arrays is the following:
\begin{equation} \label{scalingtr}
\theta_i\mapsto t\,\theta_i+c,\qquad \theta^*_i\mapsto t^*\theta^*_i+c^*,
\qquad \varphi_i\mapsto t\,t^*\varphi_i, \qquad \phi_i\mapsto t\,t^*\phi_i.
\end{equation}
Using affine transformations we can normalize a parameter array into a
convenient form. We use the following normalizations.
\begin{lemma} \label{normlps}
The general parameter arrays in {\rm\cite[Examples 35.2--35.13]{terwgen}}
can be normalized by affine transformations $(\ref{scalingtr})$ to the
following forms:
\begin{itemize}
\item The $q$-Racah case:
$\displaystyle\theta_i=s\,q^{d-2i}+\frac{q^{2i-d}}{s},\quad
\theta^*_i=s^*q^{d-2i}+\frac{q^{2i-d}}{s^*}$.
\begin{eqnarray*}
\varphi_i\equal\frac{q^{2d+2-4i}}{s\,s^*r}\left(1-q^{2i}\right)\left(1-q^{2i-2d-2}\right)
\left(s\,s^*-r\,q^{2i-d-1}\right)\left(s\,s^*r-q^{2i-d-1}\right),\\
\phi_i\equal\frac{q^{2d+2-4i}}{s\,s^*r}\left(1-q^{2i}\right)\nonumber
\left(1-q^{2i-2d-2}\right)\left(s^*r-s\,q^{2i-d-1}\right)\left(s^*-s\,r\,q^{2i-d-1}\right).
\end{eqnarray*}
\item The $q$-Hahn case: $\displaystyle\theta_i=r\,q^{d-2i},\quad
\theta^*_i=s^*q^{d-2i}+\frac{q^{2i-d}}{s^*}$,
\begin{eqnarray*}
\varphi_i\equal \frac{q^{2d+2-4i}}{r}\,
\left(1-q^{2i}\right)\left(1-q^{2i-2d-2}\right)\left(s^*r^2-q^{2i-d-1}\right),\\
\phi_i\equal-\frac{q^{d+1-2i}}{r\,s^*}\,
\left(1-q^{2i}\right)\left(1-q^{2i-2d-2}\right)\left(s^*-r^2q^{2i-d-1}\right).
\end{eqnarray*}
\item The dual $q$-Hahn case: $\displaystyle
\theta_i=s\,q^{d-2i}+\frac{q^{2i-d}}{s},\quad\theta^*_i=r\,q^{d-2i}$,
\begin{eqnarray*}
\varphi_i\equal \frac{q^{2d+2-4i}}{r}\,
\left(1-q^{2i}\right)\left(1-q^{2i-2d-2}\right)\left(s\,r^2-q^{2i-d-1}\right),\\
\phi_i\equal
\frac{q^{2d+2-4i}}{r\,s}\,\left(1-q^{2i}\right)\left(1-q^{2i-2d-2}\right)\left(r^2-s\,q^{2i-d-1}\right).
\end{eqnarray*}
\item The $q$-Krawtchouk: $\displaystyle\theta_i=q^{d-2i},\quad
\theta^*_i=s^*q^{d-2i}+\frac{q^{2i-d}}{s^*}$,
\begin{eqnarray*}
\varphi_i\equal s^*\,q^{2d+2-4i}\,\left(1-q^{2i}\right)\left(1-q^{2i-2d-2}\right),\\
\phi_i\equal\frac{1}{s^*}\,\left(1-q^{2i}\right)\left(1-q^{2i-2d-2}\right).
\end{eqnarray*}
\item The dual $q$-Krawtchouk: $\displaystyle
\theta_i=s\,q^{d-2i}+\frac{q^{2i-d}}{s},\quad\theta^*_i=q^{d-2i}$,
\begin{eqnarray*}
\varphi_i\equal s\;q^{2d+2-4i}\,\left(1-q^{2i}\right)\left(1-q^{2i-2d-2}\right),\\
\phi_i\equal\frac{q^{2d+2-4i}}{s}\,\left(1-q^{2i}\right)\left(1-q^{2i-2d-2}\right).
\end{eqnarray*}
\item The quantum $q$-Krawtchouk: $\theta_i=r\,q^{2i-d}$,\,
$\theta^*_i=r\,q^{d-2i}$,
\begin{eqnarray*}
\varphi_i\equal -\frac{q^{d+1-2i}}{r}\,\left(1-q^{2i}\right)\left(1-q^{2i-2d-2}\right),\\
\phi_i\equal \frac{q^{2d+2-4i}}{r}\,
\left(1-q^{2i}\right)\left(1-q^{2i-2d-2}\right)\left(r^3-q^{2i-d-1}\right).
\end{eqnarray*}
\item The affine $q$-Krawtchouk: $\theta_i=r\,q^{d-2i}$,\,
$\theta^*_i=r\,q^{d-2i}$,
\begin{eqnarray*}
\varphi_i\equal \frac{q^{2d+2-4i}}r\,\left(1-q^{2i}\right)\left(1-q^{2i-2d-2}\right)\left(r^3-q^{2i-d-1}\right),\\
\phi_i\equal-\frac{q^{d+1-2i}}r\,\left(1-q^{2i}\right)\left(1-q^{2i-2d-2}\right).
\end{eqnarray*}
\item The Racah case: $\theta_i=(i+u)(i+u+1)$,
$\theta^*_i=(i+u^*)(i+u^*+1)$,
\begin{eqnarray*}
\varphi_i\equal i\,(i-d-1)\,(i+u+u^*-v)\,(i+u+u^*+d+1+v),\\
\phi_i\equal i\,(i-d-1)\,(i-u+u^*+v)\,(i-u+u^*-d-1-v).
\end{eqnarray*}
\item The Hahn case: $\theta_i=i+v-\frac{d}2$,
$\theta^*_i=(i+u^*)(i+u^*+1)$,
\begin{eqnarray*}
\varphi_i\equal i\,(i-d-1)\,(i+u^*\!+2v),\\
\phi_i\equal -i\,(i-d-1)\,(i+u^*\!-2v).
\end{eqnarray*}
\item The dual Hahn case: $\theta_i=(i+u)(i+u+1)$,
$\theta^*_i=i+v-\frac{d}2$,
\begin{eqnarray*}
\varphi_i\equal i\,(i-d-1)\,(i+u+2v),\\
\phi_i\equal i\,(i-d-1)\,(i-u+2v-d-1).
\end{eqnarray*}
\item The Krawtchouk case: $\theta_i=i-\frac{d}2$, $\theta^*_i=i-\frac{d}2$,
\begin{eqnarray*}
\varphi_i\equal v\,i\,(i-d-1),\\
\phi_i\equal (v-1)\,i\,(i-d-1).
\end{eqnarray*}
\item The Bannai-Ito case: $\theta_i=(-1)^i\left(i+u-\frac{d}2\right)$,
$\theta^*_i=(-1)^i\left(i+u^*\!-\frac{d}2\right)$,
\begin{eqnarray*}
\varphi_i=\left\{ \begin{array}{cl} -i\left(i+u+u^*\!+v-\frac{d+1}2\right),&
\mbox{for $i$ even, $d$ even}.\\
-(i-d-1)\left(i+u+u^*\!-v-\frac{d+1}2\right), &
\mbox{for $i$ odd, $d$ even}.\\
-i\,(i-d-1), & \mbox{for $i$ even, $d$ odd}.\\
v^2-\left(i+u+u^*\!-\frac{d+1}2\right)^2, &
 \mbox{for $i$ odd, $d$ odd}.\end{array} \right.
\end{eqnarray*}
\begin{eqnarray*}
\phi_i=\left\{ \begin{array}{cl} i\left(i-u+u^*\!-v-\frac{d+1}2\right),& \mbox{for $i$ even, $d$ even}.\\
(i-d-1)\left(i-u+u^*\!+v-\frac{d+1}2\right), & \mbox{for $i$ odd, $d$ even}.\\
-i\,(i-d-1), & \mbox{for $i$ even, $d$ odd}.\\
v^2-\left(i-u+u^*\!-\frac{d+1}2\right)^2, &
 \mbox{for $i$ odd, $d$ odd}.\end{array} \right.
\end{eqnarray*}
\end{itemize}
In each case, $q,s,s^*,r$ are nonzero scalar parameters, or $u,u^*,v$ are
scalar parameters, such that $\theta_i\neq \theta_j$,
$\theta^*_i\neq\theta^*_j$ for $0\le i<j\le d$, and $\varphi_i\neq 0$,
$\phi_i\neq 0$ for $1\le i \le d$.
\end{lemma}
\proof These results are identical to the joint results of Lemma 6.1 and Lemma 7.1 in
\cite{normlpaw}. (Compared with the parameter arrays in \cite{terwgen},
one notable substitution is $q\mapsto q^2$. For example, to get
the normalized $q$-Racah parameter array from the general parameter array
in (\ref{qracah1})--(\ref{qracah9}), one may substitute
$q\mapsto q^2$, $s\mapsto 1/s^2q^{2d+2}$, $s^*\mapsto
1/s^*{}^2q^{2d+2}$, $r_1\mapsto r/ss^*\!q^{d+1}$ and adjust
$\theta_0,\theta^*_0,h,h^*$ by affine scalings.)
\qed\\

Affine transformations (\ref{translation}) act on Askey-Wilson relations as
well. They do not change the number of Leonard pairs with the same
Askey-Wilson relations. Hence it is enough to consider our problem for a set
of {\em normalized} Askey-Wilson relations. Possible normalizations are
discussed in \cite[Sections 4 and 8]{normlpaw}. Askey-Wilson relations
satisfied by at least one Leonard pair can be normalized as follows.
\begin{lemma} \label{normrules}
Let $AW(\beta,\gamma,\gamma^*,\varrho,\varrho^*,\omega,\eta,\eta^*)$ denote
a pair of Askey-Wilson relations satisfied by a Leonard pair. The relations
can be uniquely normalized by affine translation $(A,A^*)\mapsto
(A+c,\,A^*+c^*)$ as follows:
\begin{enumerate}
\item If $\beta\neq 2$, we can set $\gamma=0$, $\gamma^*=0$.%
\item If $\beta=2$, $\gamma\neq 0$, $\gamma^*\!\neq 0$, we can set $\varrho=0$, $\varrho^*=0$.%
\item If $\beta=2$, $\gamma=0$, $\gamma^*\!\neq 0$, we can set $\varrho^*=0$, $\omega=0$.%
\item If $\beta=2$, $\gamma^*\!=0$, $\gamma\neq 0$, we can set $\varrho=0$, $\omega=0$.%
\item If $\beta=2$, $\gamma=0$, $\gamma^*=0$ we can set $\eta=0$, $\eta^*=0$.%
\end{enumerate}
After the translation normalization, each of the two sequences
\begin{equation} \label{scalingpars}
(\gamma,\,\varrho,\,\eta,\,\eta^*)\quad\mbox{and}\quad
(\gamma^*,\varrho^*,\eta^*,\eta)
\end{equation}
contains a nonzero Askey-Wilson coefficient. By affine scaling
$(A,A^*)\mapsto (tA,\,t^*A^*)$ one can put the first nonzero coefficients
in both sequences to any convenient nonzero values.
\end{lemma}
\proof The normalization by affine translations follows from \cite[Lemma 4.1
and Part 3 of Theorem 8.1]{normlpaw}. Note that parts 6 and 7 of \cite[Lemma
4.1]{normlpaw} do not apply. Normalization by affine scaling follows from
\cite[Lemma 5.2 (or Lemma 6.2) and Lemma 7.2]{normlpaw}.
Expression (\ref{scalingpars}) is the same as \cite[formula (53)]{normlpaw}.
\qed\\

The Askey-Wilson relations for the parameter arrays of Lemma \ref{normlps}
are normalized according to the specifications of Lemma \ref{normrules}. The
following lemma presents those Askey-Wilson relations. The first nonzero
parameters in the two sequences (\ref{scalingpars}) are normalized as
in \cite[formula (54)]{normlpaw}, to the following values:
\begin{eqnarray} \label{scnormal}
\mbox{ }\gamma,\gamma^*: && \mbox{$2$ (if $\beta=2$)};\nonumber\\
\mbox{ }\varrho,\varrho^*: && \left\{\begin{array}{cl}
4\!-\!\beta^2, & \mbox{if } \beta\neq\pm 2,\\
1, & \mbox{if } \beta=\pm 2;\end{array}\right.\\
\mbox{ }\eta,\eta^*: &&
\left\{\begin{array}{cl}\sqrt{\beta\!+\!2}\,(\beta\!-\!2), & \mbox{if }
\eta\eta^*\neq 0 \mbox{ or } \omega=0,\\
\sqrt{\beta\!+\!2}\,(\beta\!-\!2)\,Q_{d+1}, & \mbox{if } \eta\eta^*=0 \mbox{
and } \omega\neq 0.\end{array}\right.\nonumber
\end{eqnarray}
We should identify $\sqrt{\beta+2}=q+1/q$. This normalization of
Askey-Wilson relations is not unique, and (in the $q$-cases) there may be
two alternative normalizations with different signs of $\sqrt{\beta+2}$; see
\cite[Section 9]{normlpaw}.
\begin{lemma} \label{awnormals}
As in the previous lemma, let $q,s,s^*,r$ denote nonzero scalar parameters,
and $u,u^*,v$ denote scalar parameters. We use the following notations:
\begin{eqnarray}
Q_j=q^j+q^{-j},\qquad Q_j^*=q^j-q^{-j}, \quad \mbox{for $j=1,2,\ldots$},\\
\label{normconsts}S=s+\frac{1}s, \qquad S^*=s^*+\frac{1}{s^*}, \qquad
R=r+\frac{1}{r}.
\end{eqnarray}
The Askey-Wilson relations for the parameter arrays of Lemma
$\ref{normlps}$ are: % the following:
\begin{itemize}
\item For the $q$-Racah case:
\begin{eqnarray}
&&\hspace{-24pt}AW\big(Q_2,\,%q^2+q^{-2},\,
0,\,0,\,-Q_2^*{}^2,\,-Q_2^*{}^2,
\,-Q_1^*{}^2\!\left(S\,S^*\!+Q_{d+1}R\right),\nonumber\\
&&Q_1Q_1^*{}^2\!\left( S\,R+Q_{d+1}S^*\right),\,
Q_1Q_1^*{}^2\!\left(S^*R+Q_{d+1}S\right)\big).
\end{eqnarray}
\item For the $q$-Hahn case:
\begin{eqnarray}
&&\hspace{-24pt}AW\big(Q_2,\,0,\,0,\,0,\,-Q^*_2{}^2,\,-Q_1^*{}^2\!\left(S^*r+Q_{d+1}r^{-1}\right),\nonumber\\
&&Q_1Q_1^*{}^2,\, Q_1Q_1^*{}^2\!\left(S^*r^{-1}+Q_{d+1}r\right)\big).
\end{eqnarray}
\item For the dual $q$-Hahn case:
\begin{eqnarray}
&&\hspace{-24pt}AW\big(Q_2,\,0,\,0,\,-Q^*_2{}^2,\,0,\,-Q_1^*{}^2\!\left(S\,r+Q_{d+1}r^{-1}\right),
\nonumber\\ &&\qquad
Q_1Q_1^*{}^2\!\left(S\,r^{-1}\!+Q_{d+1}r\right),\,Q_1Q_1^*{}^2\big).
\end{eqnarray}
\item For the $q$-Krawtchouk case:
\begin{eqnarray}
AW\big(Q_2,\,0,\,0,\,0,\,-Q_2^*{}^2,\,-Q_1^*{}^2S^*,\,0,\,Q_1Q_1^*{}^2Q_{d+1}\big).
\end{eqnarray}
\item For the dual $q$-Krawtchouk case:
\begin{eqnarray}
AW\big(Q_2,\,0,\,0,\,-Q_2^*{}^2,\,0,\,-Q_1^*{}^2S,\,Q_1Q_1^*{}^2Q_{d+1},\,0\big).
\end{eqnarray}
\item For the quantum $q$-Krawtchouk and affine $q$-Krawtchouk cases:
\begin{eqnarray}
AW\big(Q_2,\,0,\,0,\,0,\,0,\,-Q_1^*{}^2\!\left(r^2+Q_{d+1}r^{-1}\right),\,
Q_1Q_1^*{}^2,\,Q_1Q_1^*{}^2\big).\quad
\end{eqnarray}
\item For the Racah case:
\begin{eqnarray}
&&\hspace{-24pt}
AW\big(2,\,2,\,2,\,0,\,0,-2u^2-2u^*{}^2\!-2v^2-2(d\!+\!1)(u+u^*\!+v)-2d^2-4d,\nonumber\\
&&2\,u\left(u+d+1\right)\left(v-u^*\right)\left(v+u^*+d+1\right), \\
&&2\,u^*\!\left(u^*\!+d+1\right)\left(v-u\right)\left(v+u+d+1\right)\big).\nonumber
\end{eqnarray}
\item For the Hahn case:
\begin{eqnarray}
\textstyle
AW\big(2,\,0,\,2,\,1,\,0,\,0,\,-(u^*\!+1)(u^*\!+d)-2v^2-\frac{d^2}2,\,-4u^*(u^*\!+d+1)v\big).
\end{eqnarray}
\item For the dual Hahn case:
\begin{eqnarray}
\textstyle
AW\big(2,\,2,\,0,\,0,\,1,\,0,\,\,-4u(u+d+1)\,v,\,-(u+1)(u+d)-2v^2-\frac{d^2}2\big).
\end{eqnarray}
\item For the Krawtchouk case:
\begin{equation}
AW\big(2,\,0,\,0,\,1,\,1,\,2v-1,\,0,\,0\big).
\end{equation}
\item For the Bannai-Ito case, if $d$ is even:
\begin{equation}
AW\big(-2,\,0,\,0,\,1,\,1,\,4uu^*\!-2(d\!+\!1)\,v,\,2uv-(d\!+\!1)\,u^*,\,2u^*v-(d\!+\!1)\,u\big).
\end{equation}
\item For the Bannai-Ito case, if $d$ is odd:
\begin{eqnarray}
&&\hspace{-24pt}\textstyle AW\big(-2,\,0,\,0,\,1,\,1,\,
-2u^2-2u^*{}^2+2v^2+\frac{(d+1)^2}2, \nonumber \\
&&\textstyle
-u^2+u^*{}^2-v^2+\frac{(d+1)^2}4,\,u^2\!-u^*{}^2\!-v^2\!+\frac{(d+1)^2}4\big).
\end{eqnarray}
\end{itemize}
\end{lemma}
\proof These results are identical to the joint results of Lemma 6.2 and Lemma 7.2 in
\cite{normlpaw}. \qed
\\

The Askey-Wilson type can be defined for Askey-Wilson relations so that type
nominations for Leonard pairs and Askey-Wilson relations are consistent; see
\cite[Section 8]{normlpaw}. The classification of Askey-Wilson relations is
largely recapitulated by the first and third columns of Table \ref{elltab}.

An important question for us is the following. If we take concrete
Askey-Wilson relations normalized according to Lemma \ref{normrules} and
formulas (\ref{scnormal}), are all Leonard pairs satisfying them
representable by parameter arrays of Lemma \ref{normlps}? The following
lemma settles this question.
\begin{lemma} \label{normscale}
\begin{enumerate}
\item Any Leonard pair satisfying normalized Askey-Wilson relations
can be represented by a normalized parameter array, %from Lemma $(\ref{normlps1})$,
except when the Askey-Wilson type is Bannai-Ito, and $d$ is odd.%
\item Suppose that $d$ is odd. Let $(B,B^*)$ denote the Leonard pair
represented by the parameter array of the Bannai-Ito type in Lemma $\ref{normlps}$.
Then the following four Leonard pairs satisfy 
normalized Askey-Wilson relations of the Bannai-Ito type:
\begin{equation} \label{bito4}
(B,B^*),\quad (-B,B^*),\quad (B,-B^*),\quad (-B,-B^*).
\end{equation}
Of these Leonard pairs, only $(B,B^*)$ can be represented by a normalized
parameter array.
\end{enumerate}
\end{lemma}
\proof These are the results of Lemmas 9.6 in \cite{normlpaw}. The crucial observation
is that the Bannai-Ito parameter array of Lemma \ref{normlps} has the
even-indexed $\theta_i$'s and the even indexed $\theta^*_i$'s in the increasing order;
when $d$ is odd, the relations operations $\Downarrow$, $\downarrow$ preserve
this property, while affine scalings by $-1$ reverse it. \qed

\section{Correctness of Table \ref{elltab}}
\label{mainproof}

Recall that we assume $d\ge 3$. By part 1 of \cite[Theorem 8.1]{normlpaw},
all Leonard pairs satisfy Askey-Wilson relations of their own Askey-Wilson
type. Therefore, we prove correctness of Table \ref{elltab} by considering
Askey-Wilson relations of different types separately; in each case we look
only for Leonard pairs of the same Askey-Wilson type.

As mentioned just before Lemma \ref{normrules}, it is enough to consider
only normalized Askey-Wilson relations. By Lemma \ref{normscale}, all
Leonard pairs satisfying normalized Askey-Wilson relations are representable
by parameter arrays of Lemma \ref{normlps}, except when the Askey-Wilson
type is Bannai-Ito and $d$ is odd. In all cases except the Bannai-Ito case with odd $d$,
each Leonard pair solution of normalized Askey-Wilson relations is representable
by a normalized parameter array. In these cases, we just assume
free values of non-normalized coefficients in the Askey-Wilson relations of
Lemma \ref{awnormals}, equate those free values to the coefficient
expressions in the free parameters (such as $s,s^*,r$ or $u,u^*,v$) of 
the corresponding general parameter array, and count solutions of
obtained algebraic equations. We should take care of the fact  that representation
of normalized Leonard pairs by normalized parameter arrays is usually not unique.

If $\beta\neq\pm 2$, we have 4 possibilities for $q$. They are related by
the substitutions $q\mapsto-q$, $q\mapsto 1/q$ and $q\mapsto-1/q$. We may
consider $q$ fixed, because Tables 3 and 4 in \cite{normlpaw} show the
following. If a Leonard pair is represented by a $q$-parameter array of
Lemma \ref{normlps}, then it can be represented by a parameter array of
Lemma \ref{normlps} with $q$ replaced by $1/q$ as well, and such a replacement
always yields an isomorphic Leonard pair. In the $q$-Racah and, for even $d$, the
$q$-Krawtchouk and dual $q$-Krawtchouk cases, the same holds for the substitution $q\mapsto-q$. 
In the other $q$-cases, the substitution $q\mapsto -q$ leads to alternatively
normalized Askey-Wilson relations (with the other sign of $\sqrt{\beta+2}$).
In any case, it is enough to count parameter arrays for one $q$-possibility.

Other transformations of normalized parameter arrays that preserve Leonard pairs
are substitutions of their free parameters that leave the parameter arrays invariant, or
realize the $\Downarrow$-$\downarrow$-relation operators. Discarding the substitutions
which change $q$, these transformations are given in Table \ref{invartr}. 
The algebraic equations in the free parameters should be rewritten in
invariants of these transformations. Examples of these invariants (for appropriate cases)
are the expressions $S$, $S^*$, $R$ in (\ref{normconsts}).
\begin{table}
\begin{center} \begin{tabular}{|c|c|c|c|}
\hline Askey-Wilson & Parameter array & \multicolumn{2}{c|}{Conversion to relatives} \\
\cline{3-4} type & stays invariant & $\Downarrow$ &  $\downarrow$ \\ \hline %
$q$-Racah & $r\mapsto 1/r$ & $s\mapsto 1/s$ & $s^*\mapsto 1/s^*$ \\
$q$-Hahn & ---  & --- & $s^*\mapsto 1/s^*$ \\
Dual $q$-Hahn & --- & $s\mapsto 1/s$ & --- \\
$q$-Krawtchouk & ---  & --- & $s^*\mapsto 1/s^*$ \\
Dual $q$-Krawtchouk & --- & $s\mapsto 1/s$ & --- \\
Racah & $v\mapsto -v-d-1$ & $u\mapsto -u-d-1$ & $u^*\mapsto -u^*\!-d-1$ \\
Hahn & --- & --- & $u^*\mapsto -u^*\!-d-1$ \\
Dual Hahn & --- & $u\mapsto -u-d-1$ & --- \\
Bannai-Ito, $d$ odd & $v\mapsto -v$ & $u\mapsto-u$ & $u^*\mapsto -u^*$ \\
\hline
\end{tabular} \end{center}
\caption{Reparametrizations preserving Leonard pairs} \label{invartr}
\end{table}

In each Askey-Wilson case we ought to check whether solutions are generally
non-degenerate. For this, one can check generic irreducibility (over the
ring generated by free parameters) of the equation systems, or check that
degenerate solutions form subvarieties with lower dimension. For fixed
$\beta\neq\pm2$, generically degenerate Leonard pairs occur only if $\beta=2\cos\pi/j$
for some $j\in\{1,2,\ldots,d\}$, so that we have $q^{2j}=1$.

From here we consider all normalized Askey-Wilson relations case by case. We use
the notation of Lemma \ref{awnormals}. Also, we denote
\begin{equation} \label{biguv} \textstyle
U=\left(u+\frac{d+1}2\right)^2,\qquad
U^*=\left(u^*+\frac{d+1}2\right)^2,\qquad V=\left(v+\frac{d+1}2\right)^2.
\end{equation}

In the $q$-Racah case, we introduce the following indeterminants:
\begin{eqnarray} \label{qracahxyz}
x=\frac{S}{Q_{d+1}},\qquad y=\frac{S^*}{Q_{d+1}}, \qquad z=\frac{R}{Q_{d+1}}.
\end{eqnarray}
They are invariant under the relevant transformations of Table
\ref{invartr}. Equating the non-normalized Askey-Wilson coefficients gives the
equations
\begin{eqnarray} \label{qracaheq}
xy+z\equal C_1,\nonumber\\
xz+y\equal C_2,\\
yz+x\equal C_3,\nonumber
\end{eqnarray}
where
\[
C_1=-\frac{\omega}{Q^*_1{}^2\,Q_{d+1}^2},\quad
C_2=\frac{\eta}{Q_1Q^*_1{}^2\,Q_{d+1}^2},\quad
C_3=\frac{\eta^*}{Q_1Q^*_1{}^2\,Q_{d+1}^2}.
\]
Elimination of $y,z$ from (\ref{qracaheq}) gives the degree 5 equation
\begin{equation} \label{qracaheq5}
(x-C_3)(x^2-1)^2+C_1C_2\,(x^2-1)-(C_1^2+C_2^2)\,x=0.
\end{equation}
Each solution gives exactly one Leonard pair satisfying the normalized Askey-Wilson relations
$AW\left(q+q^{-1},0,0,-Q^*_2{}^2,-Q^*_2{}^2,\omega,\eta,\eta^*\right)$. There are more
solutions in terms of $(s,s^*,r)$, but distinct Leonard pairs come from
distinct $(x,y,z)$. The polynomial in (\ref{qracaheq5}) %in $x$
does not have multiple roots (in $x$) in general. Hence the generic number of Leonard pairs is 5.

In the $q$-Hahn case, invariant variables are $S^*$, % $y$ from (\ref{qracahxyz})
$r$, and free Askey-Wilson coefficients are $\omega$, $\eta^*$. Elimination
of $S^*$ gives a polynomial of degree 4 in $r$, without multiple roots in
general. The generic number of Leonard pairs is 4. The dual $q$-Hahn case is
similar.

In the $q$-Krawtchouk case, we have the equation $\omega=-Q^*_1{}^2S^*$ which
obviously has one solution in $S^*$. The dual $q$-Krawtchouk case is
similar.

For Askey-Wilson relations of the quantum/affine $q$-Krawtchouk case, we have
a cubic equation in $r$. The solutions represent 3 Leonard pairs of the same type. 
The Leonard pairs can be represented by parameter arrays of 
the quantum $q$-Krawtchouk type, or the affine $q$-Krawtchouk type.

In the Racah case, we use (\ref{biguv}) and rewrite the equations as
\begin{eqnarray} \label{racaheqs}
\omega\equal\textstyle -2\,U-2\,U^*-2\,V-\frac{(d-1)(d+3)}2,\nonumber\\
\eta\equal\textstyle 2\left(U-\frac{(d+1)^2}4\right)\left(V-U^*\right),\\
\eta^*\!\equal\textstyle2\left(U^*\!-\frac{(d+1)^2}4\right)\left(V-U\right).\nonumber
\end{eqnarray}
Here $U$, $U^*$, $V$ are invariants by Table \ref{invartr}. The degree of equations suggests
the generic number $4=1\cdot2\cdot2$ of solutions. Elimination of two of the three invariants
confirms this generic number.

In the Hahn case, we have the equations
\begin{eqnarray}
\eta\equal\textstyle -U^*-2v^2-\frac{d^2+2d-1}{4},\\
\eta^*\!\equal\textstyle -4v\left(U^*\!-\frac{(d+1)^2}4\right).\nonumber
\end{eqnarray}
The invariants are $U^*$ and $v$. Elimination of $U^*$ gives a cubic
equation in $v$:
\begin{equation} \textstyle
v^3+\left( \frac{\eta}2+\frac{d(d+2)}4 \right)v-\eta^*=0.
\end{equation}
Hence there are 3 Leonard pairs in general. The dual Hahn case is similar.

In the Krawtchouk case, we obviously have one solution.

In the Bannai-Ito case for even $d$, after setting
\begin{equation} \label{bitovars}
x=-\frac{2u}{d+1},\qquad y=-\frac{2u^*}{d+1},\qquad z=-\frac{2v}{d+1},
\end{equation}
we arrive at an equation system of the same form as in (\ref{qracaheq}), so
the generic number of solutions is 5 as well.

In the Bannai-Ito case for odd $d$ we have to keep in mind part 2 of Lemma
\ref{normscale}. The invariants under relevant transformations of Table \ref{invartr}
are $u^2$, $u^*{}^2$, $v^2$. The expressions for $\omega$, $\eta$, $\eta^*$
in Lemma \ref{awnormals} are linear in these invariants, so there is only one
solution representable by a parameter array of Lemma \ref{normlps}.
But part 2 of Lemma \ref{normscale} asserts that there are 4 Leonard pairs in total.

Correctness of Table \ref{elltab} is proved.

\section{More examples}
\label{morexamples}

First we reconsider the example in Section \ref{almostb}. 
The Askey-Wilson relations
in (\ref{bipartite1}) have the $q$-Racah type, so looking for normalized
Leonard pairs satisfying them leads us to
the equation system (\ref{qracaheq}) with $\omega=0$,
$\eta=0$, $\eta^*=0$. The equation system has the following solutions:
\begin{equation} \label{solxyz}
(x,y,z)\in\{ (0,0,0), (1,1,-1), (1,-1,1), (-1,1,1), (-1,-1,-1) \}.
\end{equation}
These solutions correspond to the Leonard pairs in (\ref{bipsols}),
respectively. Indeed, the Leonard pairs $(A_1,A_1^*)$ and $(A_2,A_2^*)$ can be obtained
from the $q$-Racah parameter array of Lemma \ref{normlps} by specializing, respectively,
\begin{equation}
s=s^*=r=\sqrt{-1}, \qquad\mbox{and}\qquad s=s^*=q^{-d-1}, \quad r=-q^{-d-1}.
\end{equation}
With this identification, the Leonard pair $(A_1,A_1^*)$ is in the [d*0*0d] basis in the terminology
of \cite{Terw24}, while the Leonard pair $(A_2,A_2^*)$ is in the [0*d*0d] basis.

\begin{example} \rm Consider the Askey-Wilson relations
$AW(-2,0,0,1,1,0,0,0)$. Leonard pairs satisfying it have the Bannai-Ito
type.

For even $d$, we have the same solutions in terms of (\ref{bitovars}) as in
(\ref{solxyz}). The solution $(0,0,0)$ corresponds to the Leonard pair
$(B_1,B^*_1)$ defined by the following matrices. The matrix for $B_1$ is
diagonal, with the following sequence of diagonal entries:
\begin{equation} \textstyle
-\frac{d}2,\quad \frac{d}2-1,\quad 2-\frac{d}2,\quad
\frac{d}2-3,\quad\ldots,\quad 1-\frac{d}2,\quad \frac{d}2.
\end{equation}
The matrix for $B_1^*$ is tridiagonal:
\begin{equation}
F_1^*= \left( \begin{array}{cccccccc} 0 & \frac{d}2 \\
\frac12 & 0 & \frac{d-1}2 \\ & 1 & 0 & \frac{d-2}2 \\
&& \frac32 & \ddots & \ddots \\
&&& \ddots & 0 & \frac12 \\ &&&& \frac{d}2 & 0 \end{array} \right).
\end{equation}
This matrix looks familiar from representation theory of the Lie algebra
$sl_2$. Let $B_1^\wedge$ denote the diagonal matrix with the same set of
diagonal entries as $B_1$, but arranged in the increasing order. Then one
can check that $(B_1^\wedge,B_1^*)$ is a Leonard pair of the Krawtchouk
type. Up to scaling, this Leonard pair (with any $d$) occurs in \cite{Go}.

For any $d$, let $(B_2,B^*_2)$ denote the Leonard pair defined by the
following matrices. The matrix for $B_2$ is diagonal, with the following
sequence of diagonal entries:
\begin{equation} \textstyle
\frac{1}2,\quad -\frac32,\quad \frac52,\quad -\frac72,\quad\ldots,\quad
(-1)^d %\varepsilon\!
\left(d\!+\!\frac12\right).
\end{equation}
%where $\varepsilon=1$.
The matrix for $B_2^*$ is tridiagonal, with exactly one nonzero entry on
the main diagonal:
\begin{equation}
\left( \begin{array}{cccccccc} \frac{d+1}2 & \frac{d}2 \\
\frac{d+2}2 & 0 & \frac{d-1}2 \\ & \frac{d+3}2 & 0 & \ddots \\
&& \ddots & \ddots & \frac32 \\ &&& d\!-\!\frac12 & 0 & 1 \\
&&&& d & 0 & \frac12 \\ &&&&& d\!+\!\frac12 & 0 \end{array} \right).
\end{equation}
It turns out that $(B_2,B^*_2)$, $(B_2,-B_2^*)$, $(-B_2,B_2^*)$,
$(-B_2,-B_2^*)$ satisfy the Askey-Wilson relations under consideration. For
even $d$, these Leonard pairs correspond to the other 4 solutions in
(\ref{solxyz}). To see the sign-flipping relation between corresponding
parameter arrays of Lemma \ref{normlps}, one has to apply the
$\Downarrow$-$\downarrow$-operations. For odd $d$, the solution
representable by the parameter array in Lemma \ref{normlps} is
$(B_2,B_2^*)$, and then we should take into account part 2 of Lemma
\ref{normscale}.
\end{example}

\begin{example} \rm
Here we consider Askey-Wilson relations of the Racah type with $\omega=0$,
$\eta=0$, $\eta^*=0$. According to (\ref{racaheqs}), there are 4 normalized
Leonard pairs satisfying these relations:
\begin{eqnarray}
(U,\,U^*,\,V)&\in&\textstyle\!\!\left\{
\left(\frac{(d+1)^2}4,\,\frac{(d+1)^2}4,\,\frac{1-6d-3d^2}4 \right),\right. \nonumber\\
&&\textstyle\left(\frac{(d+1)^2}4,\,\frac{1-6d-3d^2}4,\,\frac{(d+1)^2}4\right), \nonumber\\
&&\textstyle\left(\frac{1-6d-3d^2}4,\,\frac{(d+1)^2}4,\,\frac{(d+1)^2}4\right), \nonumber\\
&&\textstyle\left.\left(-\frac{(d-1)(d+3)}{12},\,-\frac{(d-1)(d+3)}{12},\,-\frac{(d-1)(d+3)}{12}\right)\right\}.
\end{eqnarray}
Explicit diagonal-tridiagonal forms can be obtained as follows. (They are
not necessarily normalized to standard diagonal-tridiagonal forms.) Let
$F_1$ be the diagonal matrix with the following diagonal entries:
\begin{eqnarray}
0,\quad 2,\quad 6,\quad 12,\quad 20,\quad, \ldots,\quad d(d+1).
\end{eqnarray}
Let $F_2$ denote the tridiagonal matrix with the following entries on the
superdiagonal, the main diagonal and the subdiagonal, respectively:
\begin{eqnarray}
\textstyle\frac{j-d-1}{2\,(2j-1)}\left( %\left(j+\frac{d+1}2\right)^2+\frac{3d^2+6d-1}4
j^2+(d+1)j+d^2+2d\right),\quad j=1,\ldots,d; \\
\textstyle \frac{d\,(d+1)(d+2)}2,\quad -1,\quad -3,\quad -6,\quad -10,\quad
\ldots,\quad -\frac{d(d+1)}2;\\
\textstyle\frac{j+d+1}{2\,(2j+1)}\left( %\left(j-\frac{d+1}2}\right)^2+\frac{3d^2+6d-1}4
j^2-(d+1)j+d^2+2d\right),\quad j=1,\ldots,d.
\end{eqnarray}
Let $F_2^*$ denote the matrix with the same entries as $F_2$, except that
the entry in the upper-left corner is multiplied by $-1$. Let
$\widetilde{u}=-\frac{d+1}2+\frac12\sqrt{-\frac{(d-1)(d+3)}{3}}$. Let $F_3$ be the
diagonal matrix with the diagonal entries $(j+\widetilde{u})(j+\widetilde{u}+1)$, $j=0,\ldots,1$,
and let $F_4$ be the tridiagonal matrix with the following entries on the
superdiagonal, the main diagonal and the subdiagonal, respectively:
\begin{eqnarray}
\frac{(j-d-1)\left(j+2\,\widetilde{u}\right)\left(j+3\,\widetilde{u}+d+1\right)}
{2\,\left(2j+2\,\widetilde{u}-1\right)}, \quad j=1,\ldots,d;\\ %\textstyle 
-\frac12\left(j+\widetilde{u}\right)\left(j+\widetilde{u}+1\right),\quad j=0,\ldots,d;\\
\frac{j\,\left(j+2\,\widetilde{u}+d+1\right)(j-\widetilde{u}-d-1)}
{2\,\left(2j+2\,\widetilde{u}+1\right)},\quad j=1,\ldots,d.
\end{eqnarray}
Then $(F_1,F_2)$, $(F_1,F_2^*)$, $(F_2^*,F_1)$, $(F_3,F_4)$ are matrix pairs
representing the 4 Leonard pairs.
\end{example}

\begin{example} \label{hahnexm} \rm
Here we consider Askey-Wilson relations of the Hahn type with $\eta^*=0$.
There must be solutions with $v=0$ and with $U^*\!=\frac{(d+1)^2}4$. We want
all entries in the representing matrices to be in $\QQ$, so we must have
$\left(u^*{}\!+\frac{d+1}2\right)^2\!-2v^2=\frac{(d+1)^2}{4}$. Rational
solutions of this equation can be parametrized with
$v=\frac{(d+1)t}{t^2-2}$, which gives the following family of Askey-Wilson
relations:
\begin{equation}
AW\left(2,0,2,1,0,0,\frac12-\frac{(d+1)^2(t^4+4)}{2\,(t^2-2)^2},0\right).
\end{equation}
The 3 Leonard pairs can represented by parameter arrays of Lemma
\ref{normlps} with
\begin{equation} \label{hahnsols}
(u^*,v)\in\left\{ \left(\frac{2\,(d+1)}{t^2-2},0\right),
\left(0,\frac{(d+1)\,t}{t^2-2}\right),
\left(0,-\frac{(d+1)\,t}{t^2-2}\right)\right\}.
\end{equation}

For $t=1$ we have 3 solutions, as expected. They are representable (after
the $\downarrow$ operation) by
\begin{equation}
(u^*,v)\in\left\{ \left(d+1,0\right), \left(0,d+1\right),
\left(0,-d-1\right)\right\}.
\end{equation}

For $t=3$, we have
\begin{equation}
u^*\pm 2v\in \left\{ \frac{2(d\!+\!1)}7, -\frac{9(d\!+\!1)}7,
\frac{6(d\!+\!1)}7, -\frac{13(d\!+\!1)}7, -\frac{6(d\!+\!1)}7,
-\frac{(d\!+\!1)}7\right\}.
\end{equation}
If $d+1$ is divisible by 7, then we have only one Leonard pair solution,
because two other solutions have $\varphi_i\phi_i=0$ for $i=\frac{(d+1)}7$
or $i=\frac{6(d+1)}7$ so they are degenerate. A similar statement holds for
$t=4$.

For $t=0$, all three solutions in (\ref{hahnsols}) give the Leonard pair
representable by the parameter array of Lemma \ref{normlps} with
$(u^*,v)=(0,0)$. So we have just one solution ``of multiplicity 3".
\end{example}

\begin{example} \rm Suppose that $\xi\in\CC$ satisfies $\xi^{d+1}=2$, and
consider the Askey-Wilson relations
\[ %begin{equation}
AW\left( \xi^2+\xi^{-2}, 0, 0, 0, 0, -\frac{21(\xi-\xi^{-1})^2}4,
(\xi+\xi^{-1})(\xi-\xi^{-1})^2, (\xi+\xi^{-1})(\xi-\xi^{-1})^2\right).
\] %end{equation}
Leonard pairs satisfying these relations can be represented by parameter
arrays of the quantum $q$-Krawtchouk of the affine $q$-Krawtchouk types.
There are 3 such Leonard pairs. To get affine $q$-Krawtchouk parameter
arrays, one may take $q=\xi$ so that $Q_{d+1}=\frac52$. The cubic equation
is then $r^3+\frac5{2}=\frac{21}4r$. The solutions have
$r\in\{2,\frac12,-\frac52\}$. 
\end{example}

\bibliographystyle{alpha}
\bibliography{../terwilliger,../../hypergeometric}

\begin{thebibliography}{Cur04}

\bibitem[Cur01]{Cur2hbipT}
B.~Curtin.
\newblock The {T}erwilliger algebra of a 2-homogeneous bipartite distance
  regular graph.
\newblock {\em J. Combin. Theory (B)}, 81:125--141, 2001.

\bibitem[Cur04]{Curtinpr}
B.~Curtin.
\newblock Private communication.
\newblock 2004.

\bibitem[Go02]{Go}
J.~Go.
\newblock The {T}erwilliger algebra of the {H}ypercube ${Q_D}$.
\newblock {\em Europ. J. Combinatorics}, 23:399--429, 2002.

\bibitem[KS94]{koekswart}
R.~Koekoek and R.F. Swarttouw.
\newblock The {A}skey-scheme of hypergeometric orthogonal polynomials and its
  q-analogue.
\newblock Technical Report 94-05/98-17, Delft University of Technology, {\sf
  http://aw.twi.tudelft.nl/$\sim$koekoek/askey}, 1994.

\bibitem[Ter02]{Terw24}
P.~Terwilliger.
\newblock {L}eonard pairs from 24 points of view.
\newblock {\em Rocky Mountain J. Math.}, 32(2):827--888, 2002.
\newblock {\sf http://arxiv.org/math.RA/0406577}.

\bibitem[Ter04]{TerwRacah}
P.~Terwilliger.
\newblock {L}eonard pairs and $q$-racah polynomials.
\newblock {\em Linear Algebra Appl.}, 387:235--276, 2004.
\newblock {\sf http://arxiv.org/math.RA/0306301}.

\bibitem[Ter05]{TerwLTparr}
P.~Terwilliger.
\newblock Two linear transformations each tridiagonal with respect to an
  eigenbasis of the other; comments on the parameter array.
\newblock {\em Des. Codes Cryptogr.}, 34:307--332, 2005.
\newblock {\sf http://arxiv.org/math.RA/0306291}.

\bibitem[Ter06]{terwgen}
P.~Terwilliger.
\newblock An algebraic approach to the {A}skey scheme of orthogonal
  polynomials.
\newblock In F.~Marcellan and W.~Van Assche, editors, {\em Orthogonal
  Polynomials and Special Functions: Computation and Applications}, volume 1883
  of {\em Lecture Notes in Mathematics}, pages 225--330. Springer, 2006.
\newblock {\sf http://arxiv.org/math.QA/0408390}.

\bibitem[TV04]{TerwVid}
P.~Terwilliger and R.~Vidunas.
\newblock Leonard pairs and {A}skey-{W}ilson relations.
\newblock {\em Journal of algebra and its applications}, 3(4):411--426, 2004.
\newblock {\sf http://arxiv.org/math.QA/0305356}.

\bibitem[Vid06]{normlpaw}
R.~Vidunas.
\newblock Normalized {L}eonard pairs and {A}skey-{W}ilson relations.
\newblock {\em Linear Algebra Appl.}, 2006.
\newblock {\sf http://arxiv.org/math.RA/0505041}.

\bibitem[Zhe91]{Zhidd}
A.~S. Zhedanov.
\newblock ``{H}idden symmetry'' of {A}skey-{W}ilson polynomials.
\newblock {\em Teoret. Mat. Fiz.}, 89(2):190--204, 1991.

\end{thebibliography}
%\bibliography{../hypergeometric}

\end{document}